\documentclass[11pt]{article}
\usepackage{eurosym}
\usepackage{amsmath}
\usepackage{amsfonts}
\usepackage{amssymb}
\usepackage{graphicx}
\usepackage{color}

\numberwithin{equation}{section}
\definecolor{mycolorred}{rgb}{1, 0, 0}

\newtheorem{theorem}{Theorem}[section]

\newtheorem{corollary}[theorem]{Corollary}

\newtheorem{lemma}[theorem]{Lemma}

\newtheorem{proposition}[theorem]{Proposition}
\newtheorem{remark}[theorem]{Remark}

\def\<{\langle}
\def\>{\rangle}

\begin{document}

\title{Regularization lemmas and convergence in total variation}
\author{\textsc{Vlad Bally}\thanks{
Universit\'e Paris-Est, LAMA (UMR CNRS, UPEMLV, UPEC), MathRisk INRIA,
F-77454 Marne-la-Vall\'ee, France. Email: \texttt{bally@univ-mlv.fr} }
\smallskip \\
\textsc{Lucia Caramellino}\thanks{
Dipartimento di Matematica, Universit\`a di Roma ``Tor Vergata'', and
INdAM-GNAMPA, Via della Ricerca Scientifica 1, I-00133 Roma, Italy. Email: 
\texttt{caramell@mat.uniroma2.it}.}\smallskip\\
\textsc{Guillaume Poly}\thanks{
IRMAR, Universit\'e de Rennes 1, 263 avenue du G\'en\'eral Leclerc, CS 74205
35042 Rennes, France. Email: \texttt{guillaume.poly@univ-rennes1.fr}~~~
Guillaume Poly is supported by the ANR grant UNIRANDOM.}}
\date{}
\maketitle

\begin{abstract}
We provide a simple abstract formalism of integration by parts under which
we obtain some regularization lemmas. These lemmas apply to any sequence of
random variables $(F_n)$ which are smooth and non-degenerated in some sense
and enable one to upgrade the distance of convergence from smooth
Wasserstein distances to total variation in a quantitative way. This is a
well studied topic and one can consult for instance \cite%
{[BCDi],[BKZ],[HLN], [NPy]} and the references therein for an overview of
this issue. Each of the aforementioned references share the fact that some
non-degeneracy is required along the whole sequence. We provide here the
first result removing this costly assumption as we require only
non-degeneracy at the limit. The price to pay is to control the smooth
Wasserstein distance between the Malliavin matrix of the sequence and the
Malliavin matrix of the limit which is particularly easy in the context of
Gaussian limit as their Malliavin matrix is deterministic. We then recover,
in a slightly weaker form, the main findings of \cite{[NPS]}. Another
application concerns the approximation of the semi-group of a diffusion
process by the Euler scheme in a quantitative way and under the H\"{o}%
rmander condition.
\end{abstract}

\tableofcontents

\section{Introduction}

The main historical application of Malliavin calculus, introduced in 1975 by
Paul Malliavin, was a probabilistic proof of the H\"{o}rmander regularity
criterion. But in the 40 last years it gave rise to a huge amount of various
applications, and in particular it has been developed as a branch of
stochastic analysis on the Wiener space, see the classical book of Nualart 
\cite{bib:[N]}, as well as the more recent area of research pioneered by
Nourdin and Peccati, see \cite{[NP]}. There is a major philosophical
difference between the two aforementioned views of Malliavin calculus, as
the so-called Malliavin-Stein's method (of Nourdin and Peccati), which has
been intensively studied in a recent past, mixes the formalism of
integration by parts provided by Malliavin calculus operators together with
the Stein's method. Let us recall that the quintessence of Stein's method
consists of identifying a suitable functional operator which characterizes a
specific target and use it to prove convergence towards this target in a
quantitative way. The most emblematic example is certainly the univariate
standard Gaussian distribution $\gamma$ which is characterized by the
equation $\langle f^{\prime }-x f, \gamma \rangle=0$ for every test function 
$f$. The link between Malliavin calculus appears in the identity: 
\begin{equation*}
\mathbb{E}\left( f^{\prime }(X)-X f(X)\right)=\mathbb{E}\left(f^{\prime
}(X)\left(1- \Gamma[X,-\mathcal{L}^{-1}X]\right)\right) 
\end{equation*}
where $\Gamma$ is the square field operator on the Wiener space and $%
\mathcal{L}^{-1}$ is the pseudo-inverse of the Ornstein-Uhlenbeck operator.
The quantity of interest is then $\Gamma[X,-\mathcal{L}^{-1}X]$ which is
different from the quantity $\Gamma[X,X]$ which is standard in Malliavin
calculus. Hence, although these two points of view are rather close as they
both employ the Malliavin calculus to compute distances between
distributions, they go towards different directions. Malliavin-Stein methods
focus on specific targets with specific operators in order to provide rates
of convergence whereas regularization lemmas focus on smoothness of
distribution and upgrading distances of convergence. The present article
explores this direction, namely we do not aim at proving limit theorems but
instead of that, given a limit theorem, we explore the strongest
probabilistic distances and the smoothness of the laws. In some sense, both
approaches are complementary.

\medskip

To do so, we introduce an abstract framework built on Dirichlet form theory
in which such properties may be obtained by using some integration by parts
techniques. Those techniques are very similar to the standard Malliavin
calculus but are presented in a more general framework which goes far beyond
the sole case of the Wiener space. In particular, we aim at providing a
minimalist setting leading to our regularization lemma. Our unified
framework includes the standard Malliavin calculus and different known
versions - the ``lent particle" approach for Poisson point measures
(developed by Bouleau and Denis \cite{[BD]}), the calculus based on the
splitting method developed and used in \cite{[BCDi],[CLT], [BCP]} as well as
the $\Gamma$ calculus in \cite{[BGL]}.

\medskip

The first aim of this paper is to present, in this unified framework, the
following regularization lemma (see Theorem \ref{RL}):%
\begin{equation}
\left\vert {\mathbb{E}}(f(F))-{\mathbb{E}}(f_{\delta}(F))\right\vert \leq
C\left\Vert f\right\Vert _{\infty}\Big({\mathbb{P}}(\det\sigma_{F}\leq\eta)+%
\frac{\delta^{q}}{\eta^{2q}}\mathcal{C}_{q}(F)\Big).   \label{I1}
\end{equation}
Here $f_{\delta}=f\ast\phi_{\delta}$ is the regularization by convolution by
means of a super kernel $\phi_{\delta}$ (see (\ref{kk1}) - (\ref{kk2}) and (%
\ref{super})). We use Malliavin calculus (abstract version) for $F:$ then $%
\sigma_{F}$ is the Malliavin covariance matrix and $\mathcal{C}_{q}(F)$ is a
quantity which involves the Malliavin-Sobolev norms up to order $q$ of $F.$
This inequality holds for every $\delta>0,\eta>0$ and every $q\in {\mathbb{N}%
}.$ So one may play on these parameters according to the problem at hand. A
more powerful variant of the above lemma involves derivatives of the test
function $f:$%
\begin{equation*}
\left\vert {\mathbb{E}}((\partial_{\gamma}f)(F)-{\mathbb{E}}%
((\partial_{\gamma}f_{\delta })(F)\right\vert \leq C\Big(\left\Vert
f\right\Vert _{m,\infty}{\mathbb{P}}(\det\sigma _{F}\leq\eta)+\frac{%
\delta^{q}}{\eta^{2(q+m)}}\left\Vert f\right\Vert _{\infty}\mathcal{C}%
_{q+m}(F)\Big). 
\end{equation*}
Here $\left\Vert f\right\Vert _{m,\infty}=\sum_{\left\vert \beta\right\vert
=m}\left\Vert \partial_{\beta}f\right\Vert _{\infty}$ and $m=\left\vert
\gamma\right\vert .$ Such an inequality allows to handle convergence in
distribution norms for the law of $F_{n}$ to the law of $F.$ Applications of
such convergence results are given in $\cite{[CLT],[BCP]}.$

\medskip

One important application of the regularization lemma consists in proving
that, if a sequence $F_{n}\rightarrow F$ in a distance involving smooth test
functions (as for example for the Wasserstein distance) then it converges
also in total variation distance. Of course, in order to get such a result,
we need $F_{n}$ to be smooth, in order to control $\mathcal{C}_{q}(F_{n}),$
and (more or less) non degenerated, in order to control ${\mathbb{P}}%
(\det\sigma_{F_{n}}\leq\eta).$ Actually, according to the non degeneracy
properties, several variants of the convergence result are obtained.

\medskip

Let us give an informal version of these results. Assume first that we have
the uniform non degeneracy condition $Q_{p}:=\sup_{n}{\mathbb{E}}%
((\det\sigma_{F_{n}})^{-p})<\infty$ for every $p.$ Then we prove (see Lemma %
\ref{l2} for a precise statement) that, for every given $\varepsilon>0$%
\begin{equation}
d_{TV}(F,F_{n})\leq Cd_{W}^{1-\varepsilon}(F,F_{n})   \label{I2}
\end{equation}
where $d_{TV}$ is the total variation distance and $d_{W}$ is the
Wasserstein distance. Here $C$ is a constant which depends on the Sobolev
norms and on the ``non degeneracy" constant $Q_{p}$ for some $p$ large
enough. Notice that we loose something, because we get the power $%
1-\varepsilon$ instead of $1$ for $d_{W}(F,F_{n}).$ This is somehow a
technical drawback of our method which is based on an optimization
procedure. A more careful examination of this optimization procedure is
likely to provide logarithmic losses but this would result in highly
technical computations which fall beyond the scope of this paper. Let us
emphasize that the previous estimate requires non-degeneracy assumptions
along the whole sequence $(F_n)$ which may be in general rather hard to
check. Assumptions of non-degeneracy on the sequence $(F_n)$ may sometimes
be provided by classic anti-concentration estimates. For instance, when the
underlying Gaussian functionals are polynomials, the Carbery-Wright estimate
gives a kind of non-degeneracy but in a much weaker way. The reader can
consult \cite{[BKZ],[NPy]} for results in this direction. Another reference
of interest is \cite{[HLN]} where convergence of densities is explored when
the limit is Gaussian and under the same non-degeneracy assumption. Finally,
let us mention the reference \cite{[Poly]} which shows that in the
particular setting of quadratic forms of Gaussian vectors, the central
convergence automatically implies the required non-degeneracy assumptions
and the previous results apply under the sole assumption of Gaussian
convergence.

\medskip

In order to bypass this major issue, we are able to obtain a variant of the
above estimate without assuming nothing on the non degeneracy of $F_{n}$ (so 
$Q_{p}$ may be infinite). In Proposition \ref{l3} we prove that, for every $%
\varepsilon>0$ 
\begin{equation}
d_{TV}(F,F_{n})\leq C(d_{W}^{1-\varepsilon}(F,F_{n})+d_{W}^{1-\varepsilon
}(\det\sigma_{F},\det\sigma_{F_{n}}))   \label{I3}
\end{equation}
where $C$ depends on the Sobolev norms and $\varepsilon$ only (and not on
the non degeneracy constant $Q_{p}).$ And in Proposition \ref{l3bis} we
prove the same result with $\det\sigma_F$ replaced by a general non
degenerate positive random variable, so one just needs that $\sigma_{F_{n}}$
converges (non necessarily to $\sigma_{F}$).

In concrete examples it may be difficult to precisely estimate $%
d_{W}(\det\sigma_{F},\det\sigma_{F_{n}}),$ but then one may use the standard
upper bound $d_{W}(\det\sigma_{F},\det\sigma_{F_{n}})\leq C\left\Vert
DF-DF_{n}\right\Vert_1 .$ Doing this is not innocent, because on replace
``weak distances" with ``strong" onces and this may induce a serious loss of
accuracy: for example the weak distance is of order $\frac{1}{n}$ while the
strong distance is $\frac{1}{\sqrt{n}}.$ However, the aforementioned result
completely covers the case of central convergence as in this case $\sigma_F$
is a deterministic matrix and the quantity $d_{W}(F,F_n)$ is easy to
estimate. Using this strategy we recover a central result of Nourdin-Peccati
theory \cite{[NPS]} establishing multivariate total variation estimates for
suitable sequences converging to Gaussian. Proofs are completely different
as the proof of the aforementioned article employs tools of information
theory and provides stronger results such as convergence in entropy. On the
other hand, our result is more general and requires much less structural
information on the sequence approximating the Gaussian law.

\medskip

We also illustrate the above results in the framework of the approximation
of the semi-group of a diffusion process by using the Euler scheme: if one
assumes uniform ellipticity, then one has uniform non degeneracy for the
Euler scheme and may use (\ref{I2}). But if one works under H\"{o}rmander
condition, then the Euler scheme is degenerated so $Q_{p}=\infty.$ In \cite%
{[BT]} this problem has been discussed and the authors have been obliged to
work with a slightly regularized Euler scheme in order to bypass this
difficulty. Now, one may use (\ref{I3}) and to get the convergence for the
real Euler scheme (without regularization). But one looses accuracy: we pass
from $\frac{1}{n}$ to $\frac{1}{\sqrt{n}},$ so the result is not optimal.

A last type of results concerns the distance between density functions. This
issue has already been discussed in \cite{[BCDi]}. Here, in the Theorem \ref%
{P} we prove the following: if $F$ and $G$ are smooth and non degenerated
then the density functions $p_{F}$ and $p_{G}$ exists and are smooth.
Moreover, for every multi index $\alpha$ and for every $\varepsilon>0$%
\begin{equation}
\left\Vert \partial^{\alpha}p_{F}-\partial^{\alpha}p_{G}\right\Vert _{\infty
}\leq Cd_W(F,G)^{1-\varepsilon}.   \label{I4}
\end{equation}
This is a striking improvement with respect to the estimate obtained in \cite%
{[BCDi]}, see (2.53) there. %

\section{Abstract framework}

In this section we present an abstract framework which covers most of the
known variants of Malliavin calculus and which allows to obtain the
integration by parts formula that we need.

We consider a probability space $(\Omega ,\mathcal{F},{\mathbb{P}})$ and a
subset $\mathcal{E}\subset \cap _{p>1}L^{p}(\Omega ;{\mathbb{R}}).$ The
guiding example is $\mathcal{E}=\mathcal{S}$ or also, $\mathcal{E}={\mathbb{D%
}}^{\infty }$ (the space of simple functionals respectively the space of
smooth functionals in the classical Malliavin calculus). We assume that for
every $\phi \in C_{p}^{\infty }({\mathbb{R}}^{d})$ (sooth functions with
polynomial growth) and every $F\in \mathcal{E}^{d},$ one has $\phi (F)\in 
\mathcal{E}$. In particular $\mathcal{E}$ is an algebra. In the sequel we
will also use the following consequence. For $\eta >0$ we denote by $\Psi
_{\eta }:(0,\infty )\rightarrow {\mathbb{R}}$ a smooth function which is
equal to zero on $(0,\eta /2)$ and to one on $(\eta ,\infty ).$ Then for
every $\eta >0,$ 
\begin{equation}
F\in \mathcal{E\quad }\Longrightarrow \quad \frac{1}{F}\Psi _{\eta }(F)\in 
\mathcal{E}  \label{invers}
\end{equation}
Moreover we consider

\begin{itemize}
\item[$\clubsuit$] $\Gamma:\mathcal{E}\times\mathcal{E}\rightarrow \mathcal{E%
}$ which is a symmetric bilinear form such that $\Gamma(F,F)\geq0$ and $%
\Gamma(F,F)=0$ iff $F=0.$

In the language of Dirichlet forms $\Gamma$ is the \textit{carr\'{e} du
champs} operator. Notice that, since $\Gamma(F,G)\in\mathcal{E}$ and $%
\mathcal{E}$ is an algebra, if $F,G,H\in\mathcal{E}$ then $\Gamma(F,G)H\in%
\mathcal{E}.$ We also may define $\Gamma(F,\Gamma(G,H)).$

\item[$\clubsuit$] $L:\mathcal{E}\rightarrow \mathcal{E}$ which is a linear
operator.
\end{itemize}

We assume:

\begin{itemize}
\item[$\blacklozenge$] \textbf{[Chain rule]} For every $\phi\in
C_{p}^{\infty}({\mathbb{R}}^{d})$ and $F=(F_{1},...,F_{d})\in\mathcal{E}^d$%
\begin{equation}
\Gamma(\phi(F),G)=\sum_{i=1}^{d}\partial_{i}\phi(F)\Gamma(F_{i},G)
\label{Chain}
\end{equation}
In particular, taking $\phi(x,y)=xy$ we obtain 
\begin{equation}
\Gamma(FH,G)=F\Gamma(H,G)+H\Gamma(F,G)  \label{Prod}
\end{equation}

\item[$\blacklozenge$] \textbf{[Duality formula]} For every $F,G\in\mathcal{E%
}$, 
\begin{equation}
{\mathbb{E}}(\Gamma(F,G))=-{\mathbb{E}}(FLG)=-{\mathbb{E}}(GLF).
\label{DUAL}
\end{equation}
Notice that we also have the following extension of the duality formula:
using the duality first and the chain rule for the function $\phi(x,y)=xy$
we get for every $F,G,H\in\mathcal{E}$:%
\begin{equation*}
{\mathbb{E}}(HFLG)={\mathbb{E}}(\Gamma(HF,G))={\mathbb{E}}%
(H\Gamma(F,G)+F\Gamma(H,G))
\end{equation*}
so that 
\begin{equation}
{\mathbb{E}}(\Gamma(F,G)H)={\mathbb{E}}(F(HLG-\Gamma(H,G))  \label{dual}
\end{equation}
\end{itemize}

This gives the standard integration by parts formula that we present now.

\begin{lemma}
Let $F=(F_{1},...,F_{d})\in\mathcal{E}^{d}$ and let $\sigma_{F}^{i,j}=%
\Gamma(F_{i},F_{j}),i,j=1,...,d$ be its Malliavin covariance matrix. We
suppose that $\sigma_{F}$ is invertible, we denote $\gamma_{F}=%
\sigma_{F}^{-1},$ and we assume that 
\begin{equation}
\gamma_{F}^{k,i}\in\mathcal{E\quad}\forall i,k=1,...,d.  \label{ND}
\end{equation}
Then for every $\phi\in C_{p}^{\infty}({\mathbb{R}}^{d})$ and $G\in\mathcal{E%
}$%
\begin{equation}
{\mathbb{E}}(\partial_{i}\phi(F)G)={\mathbb{E}}(\phi(F)H_{i}(F,G)\quad with
\label{IP}
\end{equation}
with%
\begin{equation}
H_{i}(F,G)=\sum_{k=1}^{d}G(\gamma_{F}^{k,i}LF_{k}-\Gamma(%
\gamma_{F}^{k,i},F_{k}))-\sum_{k=1}^{d}\gamma_{F}^{k,i}\Gamma(G,F_{k}).
\label{W}
\end{equation}
Moreover, iterating this relation we get%
\begin{equation}
{\mathbb{E}}(\partial_{\alpha}\phi(F)G)={\mathbb{E}}(FH_{\alpha}(F,G))\quad
with  \label{Iteration}
\end{equation}
with $H_{\alpha}(F,G)$ obtained by iterations: if $\alpha=(\alpha
_{1},...,\alpha_{m})\in\{1,...,d\}^{m}$ and\ $\overline{\alpha}=(\alpha
_{1},...,\alpha_{m-1})$ then we define $H_{\alpha}(F,G)=H_{\alpha_{m}}(F,H_{%
\overline{\alpha}}(F,G)).$
\end{lemma}

\begin{remark}
If $\det\sigma_{F}$ is almost surely invertible and $(\det\sigma_{F})^{-1}\in%
\mathcal{E}$ then (\ref{ND}) is verified
\end{remark}

\textbf{Proof}. We use the chain rule and we get%
\begin{equation*}
\Gamma(\phi(F),F_{k})=\sum_{i=1}^{d}\partial_{i}\phi(F)\Gamma(F_{i},F_{k})=%
\sum_{i=1}^{d}\partial_{i}\phi(F)\sigma_{F}^{i,k}.
\end{equation*}
This gives $\nabla\phi(F)=\Gamma(\phi(F),F)\gamma_{F}$ which, on components
reads%
\begin{equation*}
\partial_{i}\phi(F)=\sum_{k=1}^{d}\Gamma(\phi(F),F_{k})\gamma_{F}^{k,i}.
\end{equation*}
Then, by (\ref{dual}) 
\begin{align*}
{\mathbb{E}}(\partial_{i}\phi(F)G) & =\sum_{k=1}^{d}{\mathbb{E}}%
(\Gamma(\phi(F),F_{k})\gamma _{F}^{k,i}G) \\
& =\sum_{k=1}^{d}{\mathbb{E}}(\phi(F)(\gamma_{F}^{k,i}GLF_{k}-\Gamma(%
\gamma_{F}^{k,i}G,F_{k})).
\end{align*}
$\square$\medskip

The non degeneracy hypothesis considered in the previous lemma is sometimes
too strong (this is the case in our framework). So we present now a
localized version of the previous integration by parts formula. We recall
that $\Psi_{\eta}(x)$ is a smooth function which is null for $x\leq\eta/2$
and equal to one for $x>\eta.$ Notice that $\sigma_{F}$ is invertible on the
set $\{\Psi_{\eta}(\det\sigma_{F})>0\}$ so we are able to define%
\begin{equation*}
\gamma_{F,\eta}^{i,j}:=\gamma_{F}^{i,j}\Psi_{\eta}(\det\sigma_{F}).
\end{equation*}
And, by (\ref{invers}) we know that $\gamma_{F,\eta}^{i,j}\in\mathcal{E}.$

\begin{lemma}
Let $F=(F_{1},...,F_{d})\in\mathcal{E}^{d}.$ Then for every $\phi\in
C_{p}^{\infty}({\mathbb{R}}^{d})$ and $G\in\mathcal{E}$%
\begin{equation}
{\mathbb{E}}(\partial_{i}\phi(F)G\Psi_{\eta}(\det\sigma_{F}))={\mathbb{E}}%
(FH_{\eta,i}(F,G))\quad with  \label{localIP}
\end{equation}
with 
\begin{equation}
H_{\eta,i}(F,G)=\sum_{k=1}^{d}G(\gamma_{F,\eta}^{k,i}LF_{k}-\Gamma
(\gamma_{F,\eta}^{k,i}G,F_{k})).  \label{LocW}
\end{equation}
Moreover, iterating this relation we get%
\begin{equation}
{\mathbb{E}}(\partial_{\alpha}\phi(F)G\Psi_{\eta}(\det\sigma_{F}))={\mathbb{E%
}}(FH_{\eta,\alpha }(F,G))\quad with  \label{iteration}
\end{equation}
with $H_{\eta,\alpha}(F,G)$ obtained by iterations: if $\alpha=(\alpha
_{1},...,\alpha_{m})\in\{1,...,d\}^{m}$ and\ $\overline{\alpha}=(\alpha
_{1},...,\alpha_{m-1})$ then we define $H_{\eta,\alpha}(F,G)=H_{\eta
,\alpha_{m}}(F,H_{\eta,\overline{\alpha}}(F,G)).$
\end{lemma}

\textbf{Proof}. The proof is almost the same as above. The only change is
that in the first step we multiply with $\Psi_{\eta}(\det\sigma_{F})$ and we
write%
\begin{equation*}
\Gamma(\phi(F),F_{k})\Psi_{\eta}(\det\sigma_{F})=\sum_{i=1}^{d}\partial
_{i}\phi(F)\sigma_{F}^{i,k}\Psi_{\eta}(\det\sigma_{F}).
\end{equation*}
On the set $\Psi_{\eta}(\det\sigma_{F})>0$ the matrix $\sigma_{F}$ is
invertible so we get 
\begin{equation*}
\partial_{i}\phi(F)\Psi_{\eta}(\det\sigma_{F})=\sum_{k=1}^{d}\Gamma
(\phi(F),F_{k})\Psi_{\eta}(\det\sigma_{F})\gamma_{F,\eta}^{k,i}=%
\sum_{k=1}^{d}\Gamma(\phi(F),F_{k})\gamma_{F,\eta}^{k,i}
\end{equation*}
and then, by (\ref{dual}) 
\begin{align*}
{\mathbb{E}}(\partial_{i}\phi(F)G\Psi_{\eta}(\det\sigma_{F})) &
=\sum_{k=1}^{d}{\mathbb{E}}(\Gamma(\phi(F),F_{k})\gamma_{F,\nu}^{k,i}G \\
& =\sum_{k=1}^{d}{\mathbb{E}}(\phi(F)(\gamma_{F,\eta}^{k,i}GLF_{k}-\Gamma(%
\gamma _{F,\eta}^{k,i}G,F_{k})).
\end{align*}
$\square$

\subsection{Norms}

In order to be able to give estimates of $H_{\alpha}(F,G)$ we need to assume
that $\Gamma$ is given by a derivative operator as follows (this is actually
always true, see Mokobodzki \cite{[M]}):

\begin{itemize}
\item[$\clubsuit$] We assume that there exists a separable Hilbert space $%
\mathcal{H}$ and a linear application $D:\mathcal{E}\rightarrow\cap
_{p>1}L^{p}(\Omega;\mathcal{H})$ such that%
\begin{align*}
i)&\qquad\Gamma(F,G) =\left\langle DF,DG\right\rangle _{\mathcal{H}} \\
ii)&\qquad D_{h}F :=\left\langle DF,h\right\rangle \in\mathcal{E}.
\end{align*}
We also assume that we have the chain rule: for $\phi\in C^{\infty}({\mathbb{%
R}}^{d})$ and $F\in\mathcal{E}^{d}$ we have%
\begin{align*}
iii)&\qquad D\phi(F)=\sum_{i=1}^{d}\partial_{i}\phi(F)DF_{i}.
\end{align*}
\end{itemize}

Then we may define higher order derivatives in the following way. $D^{2}:%
\mathcal{E}\rightarrow\cap_{p>1}L^{p}(\Omega;\mathcal{H}^{\otimes2})$ is
defined by $\left\langle D^{2}F,h_{1}\otimes h_{2}\right\rangle _{\mathcal{H}%
^{\otimes2}}=D_{h_{2}}D_{h_{1}}F.$ So, if we denote $D_{h_{1},h_{2}}^{2}F=%
\left\langle D^{2}F,h_{1}\otimes h_{2}\right\rangle _{\mathcal{H}%
^{\otimes2}} $ then $D_{h_{1},h_{2}}^{2}F=D_{h_{2}}D_{h_{1}}F.$ In a similar
way we define (by recurrence)%
\begin{equation*}
D_{h_{1},...,h_{k}}^{k}F=D_{h_{k}}D_{h_{1},...,h_{k-1}}^{k-1}F.
\end{equation*}
We introduce now the norms%
\begin{equation*}
\left\vert F\right\vert _{1,k}=\sum_{i=1}^{k}\left\vert D^{i}F\right\vert _{%
\mathcal{H}^{\otimes i}},\qquad\left\vert F\right\vert _{k}=\left\vert
F\right\vert +\left\vert F\right\vert _{1,k}=\left\vert F\right\vert
+\sum_{i=1}^{k}\left\vert D^{i}F\right\vert _{\mathcal{H}^{\otimes i}}
\end{equation*}
For $F=(F_{1},...,F_{d})\in\mathcal{E}^{d}$ we define 
\begin{equation*}
\left\vert F\right\vert _{1,k}=\sum_{i=1}^{d}\left\vert F_{i}\right\vert
_{1,k},\qquad\left\vert F\right\vert _{k}=\sum_{i=1}^{d}\left\vert
F_{i}\right\vert _{k}.
\end{equation*}

Notice that since $\mathcal{H}$ is separable we may take an orthonormal base 
$(e_{i})_{i\in {\mathbb{N}}}$ and denote $D_{i}F=D_{e_{i}}F=\left\langle
DF,e_{i}\right\rangle .$ Then $DF=\sum_{i=1}^{\infty }D_{i}F\times e_{i}$
and more generally%
\begin{equation*}
D^{k}F=\sum_{i_{1},...,i_{k}}D_{i_{1},...,i_{k}}F\times \otimes
_{j=1}^{k}e_{j}.
\end{equation*}%
Moreover we denote%
\begin{equation}
\alpha _{k}=\frac{\left\vert F\right\vert _{1,k+1}^{2(d-1)}(\left\vert
F\right\vert _{1,k+1}+\left\vert LF\right\vert _{k})}{\det \sigma _{F}}%
,\qquad \beta _{k}=\frac{\left\vert F\right\vert _{1,k+1}^{2d}}{\det \sigma
_{F}}  \label{norm1}
\end{equation}%
and%
\begin{align}
\mathcal{K}_{n,k}(F)& =(\left\vert F\right\vert _{1,k+n+1}+\left\vert
LF\right\vert _{k+n})^{n}(1+\left\vert F\right\vert _{1,k+n+1})^{2d(2n+k)},
\label{norm2} \\
\mathcal{C}_{n}(F)& =\mathcal{K}_{n,0}(F)=(\left\vert F\right\vert
_{1,n+1}+\left\vert LF\right\vert _{n})^{n}(1+\left\vert F\right\vert
_{1,n+1})^{4dn}  \label{norm2'} \\
\mathcal{C}_{n,p}(F)& =\left\Vert \mathcal{C}_{n}(F)\right\Vert _{p}.
\label{norm2''}
\end{align}%
Then, the following lemma is proved in the Appendix of \cite{[CLT]}:

\begin{lemma}
\label{L1} \textbf{A. }Let $F\in\mathcal{E}^{d}.$ Suppose that $\det\sigma
_{F}(\omega)>0$. Then, for every $k$ and $n$ there exists a constant $C$
such that for every multi-index $\rho$ with $\left\vert \rho\right\vert \leq
n$ one has 
\begin{equation}
\left\vert H_{\rho}(F,G)\right\vert _{k}\leq C\alpha_{k+n}^{n}\sum
_{p_{1}+p_{2}\leq k+n}\left\vert G\right\vert
_{p_{2}}(1+\beta_{k+n})^{p_{1}}.  \label{norm3}
\end{equation}
\textbf{B.} For every $\eta>0$ 
\begin{equation}
\left\vert H_{\rho}(F,\Psi_{\eta}(\det\sigma_{F})G)\right\vert _{k}\leq 
\frac{C}{\eta^{2n+k}}\times\mathcal{K}_{n,k}(F)\times\left\vert G\right\vert
_{k+n}.  \label{norm4}
\end{equation}
\end{lemma}

As an immediate consequence of (\ref{norm4}) and of (\ref{iteration}) we get

\begin{corollary}
Let $F\in \mathcal{E}^{d}$ and $\eta >0.$ Then%
\begin{equation}
\left\vert {\mathbb{E}}(\partial _{\alpha }f(F)\Psi _{\eta }(\det \sigma
_{F}))\right\vert \leq \frac{C}{\eta ^{2\left\vert \alpha \right\vert }}%
\times \mathcal{C}_{\left\vert \alpha \right\vert ,1}(F)\times \left\Vert
f\right\Vert _{\infty }.  \label{norm5}
\end{equation}
\end{corollary}

\bigskip

\section{Regularization Lemma}

We go now on and we give the regularization lemma. We recall that a super
kernel $\phi :{\mathbb{R}}^{d}\rightarrow {\mathbb{R}}$ is a function which
belongs to the Schwartz space $\mathbf{S}$\ (infinitely differentiable
functions wit rapid decrease) and such that for every multi-indexes $\alpha $
and $\beta $, one has 
\begin{align}
& \int \phi (x)dx=1\mbox{ and }\int y^{\alpha }\phi (y)dy=0 \mbox{ for }
|\alpha |\geq 1,  \label{kk1} \\
& \int \left\vert y\right\vert ^{m}\left\vert \partial _{\beta }\phi
(y)\right\vert dy<\infty .  \label{kk2}
\end{align}%
As usual, if $\alpha =(\alpha _{1},....,\alpha _{m})$ then $y^{\alpha
}=\prod_{i=1}^{m}y_{\alpha _{i}}$. Since super kernels play a crucial role
in our approach we give here the construction of such an object. We do it in
dimension $d=1$ and then we take tensor products. We take $\psi \in \mathbf{S%
}$ which is symmetric and equal to one in a neighborhood of zero and we
define $\phi =\mathcal{F}^{-1}\psi ,$ the inverse of the Fourier transform
of $\psi .$ Since $\mathcal{F}^{-1}$ sends $\mathbf{S}$\ into $\mathbf{S}$
the property (\ref{kk2}) is verified. And we also have $0=\psi
^{(m)}(0)=i^{-m}\int x^{m}\phi (x)dx$ so (\ref{kk1}) holds as well. We
finally normalize in order to obtain $\int \phi =1.$

We fix a super kernel $\phi$. For $\delta\in(0,1)$ and for a function $f$ we
define 
\begin{equation}  \label{super}
\phi_{\delta}(y)=\frac{1}{\delta^{d}}\phi\Big(\frac{y}{\delta}\Big)\quad %
\mbox{and}\quad f_{\delta}=f\ast\phi_{\delta},
\end{equation}
the symbol $\ast$ denoting convolution. Moreover, for $f\in C_{b}^{\infty }({%
\mathbb{R}}^{d})$ we denote%
\begin{equation}
\left\Vert f\right\Vert _{k,\infty}=\sum_{0\leq\left\vert \alpha\right\vert
\leq k}\left\Vert \partial_{\alpha}f\right\Vert _{\infty}.  \label{kk3}
\end{equation}

\begin{lemma}
\label{RL}Let $f\in C_{b}^{\infty }({\mathbb{R}}^{d})$ and $F\in \mathcal{E}%
^{d}.$ For every $q,m\in {\mathbb{N}}$\ there exists a universal constant $C 
$ (depending on $q$ and $m$ only) such that for every multi-index $\gamma $
with $\left\vert \gamma \right\vert =m,$ every $\delta >0$ and every $\eta
>0 $%
\begin{equation}
\left\vert {\mathbb{E}}(\partial _{\gamma }f(F))-{\mathbb{E}}(\partial
_{\gamma }f_{\delta }(F))\right\vert \leq C\left\Vert \partial _{\gamma
}f\right\Vert _{\infty }{\mathbb{P}}(\det \sigma _{F}\leq \eta )+\frac{%
\delta ^{q}}{\eta ^{2(q+m)}}\left\Vert f\right\Vert _{\infty }\mathcal{C}%
_{q+m,1}(F)  \label{e6}
\end{equation}%
with $\mathcal{C}_{q+m}(F)$ given in (\ref{norm2'}).\ In particular, taking $%
m=0$%
\begin{equation}
\left\vert {\mathbb{E}}(f(F))-{\mathbb{E}}(f_{\delta }(F))\right\vert \leq
C\left\Vert f\right\Vert _{\infty }({\mathbb{P}}(\det \sigma _{F}\leq \eta )+%
\frac{\delta ^{q}}{\eta ^{2q}}\mathcal{C}_{q,1}(F))  \label{e7}
\end{equation}
\end{lemma}

\begin{remark}
A similar estimate holds for $\left\vert {\mathbb{E}}(G\partial _{\gamma
}f(F))-{\mathbb{E}}(G\partial _{\gamma }f_{\delta }(F))\right\vert $ with $%
G\in \mathcal{E}.$ But in this case one has to replace ${\mathbb{P}}(\det
\sigma _{F}\leq \eta )$ with $\left\Vert G\right\Vert _{2}{\mathbb{P}}%
^{1/2}(\det \sigma _{F}\leq \eta )$\ and $\mathcal{C}_{q+m,1}(F)$ by $%
\left\Vert \left\vert G\right\vert _{q+m+1}\right\Vert _{2}\mathcal{C}%
_{q+m,2}(F)$ in the right hand side of (\ref{e6}). The proof is the same.
\end{remark}

\textbf{Proof}. The proof is given in \cite{[CLT]} in a particluar
framework, but, for the convenience of the reader, we recall it here. Using
Taylor expansion of order $q$ (with integral remainder) and (\ref{kk1}) we
obtain 
\begin{equation*}
\partial _{\gamma }f(x)-\partial _{\gamma }f_{\delta }(x)=\int (\partial
_{\gamma }f(x)-\partial _{\gamma }f(y))\phi _{\delta }(x-y)dy=\int R_{\gamma
,q}(x,y)\phi _{\delta }(x-y)dy
\end{equation*}%
with 
\begin{equation*}
R_{\gamma ,q}(x,y)=\frac{1}{q!}\sum_{\left\vert \alpha \right\vert
=q}\int_{0}^{1}\partial _{\alpha }\partial _{\gamma }f(x+\lambda
(y-x))(x-y)^{\alpha }(1-\lambda )^{q}d\lambda .
\end{equation*}%
By a change of variable we get 
\begin{equation*}
\int R_{\gamma ,q}(x,y)\phi _{\delta }(x-y)dy=\frac{1}{q!}\sum_{\left\vert
\alpha \right\vert =q}\int_{0}^{1}\int dz\phi _{\delta }(z)\partial _{\alpha
}\partial _{\gamma }f(x+\lambda z)z^{\alpha }(1-\lambda )^{q}d\lambda .
\end{equation*}%
So, we have 
\begin{align}
& {\mathbb{E}}(\Psi _{\eta }(\det \sigma _{F})\partial _{\gamma }f(F))-{%
\mathbb{E}}(\Psi _{\eta }(\det \sigma _{F})\partial _{\gamma }f_{\delta }(F))
\label{e7b} \\
& ={\mathbb{E}}\Big(\int \Psi _{\eta }(\det \sigma _{F})R_{\gamma
,q}(F,y)\phi _{\delta }(F-y)dy\Big)  \notag \\
& =\frac{1}{q!}\sum_{\left\vert \alpha \right\vert =q}\int_{0}^{1}\int
dz\phi _{\delta }(z){\mathbb{E}}\big(\Psi _{\eta }(\det \sigma _{F})\partial
_{\alpha }\partial _{\gamma }f(F+\lambda z)\big)z^{\alpha }(1-\lambda
)^{q}d\lambda .  \notag
\end{align}%
As a consequence of (\ref{norm5})%
\begin{equation*}
\left\vert {\mathbb{E}}(\Psi _{\eta }(\det \sigma _{F})\partial _{\alpha
}\partial _{\gamma }f(F+\lambda z))\right\vert \leq \frac{C}{\eta ^{2(q+m)}}%
\mathcal{C}_{q+m,1}(F)\left\Vert f\right\Vert _{\infty }.
\end{equation*}%
We also have, if $\left\vert \alpha \right\vert =q,$\ 
\begin{equation}
\int dz\left\vert \phi _{\delta }(z)z^{\alpha }\right\vert \leq \delta
^{q}\int \left\vert \phi (z)\right\vert \left\vert z\right\vert ^{|\alpha
|}dz  \label{e7e}
\end{equation}%
so we conclude that 
\begin{equation*}
\left\vert {\mathbb{E}}(\Psi _{\eta }(\det \sigma _{F}))\partial _{\gamma
}f(F))-{\mathbb{E}}(\Psi _{\eta }(\det \sigma _{F}))\partial _{\gamma
}f_{\delta }(F))\right\vert \leq \frac{C\delta ^{q}}{\eta ^{2(q+m)}}\mathcal{%
C}_{q+m,1}(F)\left\Vert f\right\Vert _{\infty }.
\end{equation*}

We write now 
\begin{equation*}
\left\vert {\mathbb{E}}((1-\Psi_{\eta}(\det\sigma_{F}))\partial_{\gamma
}f(F))\right\vert \leq\left\Vert \partial_{\gamma}f\right\Vert _{\infty}{%
\mathbb{P}}(\det\sigma_{F}<\eta).
\end{equation*}
A similar inequality holds for $\partial_{\gamma}f_{\delta}$ and so we
conclude. $\square$

\medskip

Under a weak non degeneracy condition for $F$, we get the following
immediate consequence.

\begin{corollary}
\label{C1}Let $F\in \mathcal{E}^{d}.$ Suppose that for some $\kappa >0$ one
has 
\begin{equation}
{\mathbb{P}}(\det \sigma _{F}\leq \eta )\leq \theta _{\kappa }(F)\eta
^{\kappa }\qquad \forall \eta >0.  \label{e8a}
\end{equation}%
Then for every, $q\in {\mathbb{N}}$ and $\delta >0$%
\begin{equation}
\left\vert {\mathbb{E}}(f(F))-{\mathbb{E}}(f_{\delta }(F))\right\vert \leq
C\left\Vert f\right\Vert _{\infty }\mathcal{C}_{q,1}^{\frac{\kappa }{\kappa
+2q}}(F)\theta _{\kappa }^{\frac{2q}{\kappa +2q}}(F)\times \delta ^{\frac{%
\kappa q}{\kappa +2q}}.  \label{e8b}
\end{equation}
\end{corollary}

\textbf{Proof}. One plugs (\ref{e8a}) into (\ref{e7}) and then optimize over 
$\eta$ and $\delta.$ $\square$\medskip

\textbf{Example} Let $F=\Delta ^{2}$ with $\Delta $ a standard normal random
variable. We use standard Malliavin calculus to see what (\ref{e8b}) gives
in this case. We have $DF=2\Delta $ so that $\sigma _{F}=4\Delta ^{2}$ and
then (\ref{e8a}) reads%
\begin{equation*}
{\mathbb{P}}(4\Delta ^{2}\leq \eta )={\mathbb{P}}(\left\vert \Delta
\right\vert \leq \frac{1}{2}\sqrt{\eta })\sim \eta ^{1/2}.
\end{equation*}%
So $\kappa =\frac{1}{2}$ here and (\ref{e8b}) gives (one takes $q\rightarrow
\infty )$ 
\begin{equation*}
\left\vert {\mathbb{E}}(f(F))-{\mathbb{E}}(f_{\delta }(F))\right\vert \leq
C\left\Vert f\right\Vert _{\infty }\times \delta ^{\frac{\kappa }{2}%
-}=C\left\Vert f\right\Vert _{\infty }\times \delta ^{\frac{1}{4}-}.
\end{equation*}%
But some informal computations give 
\begin{equation*}
\left\vert {\mathbb{E}}(f(F))-{\mathbb{E}}(f_{\delta }(F))\right\vert \sim
C\left\Vert f\right\Vert _{\infty }\times \delta ^{\frac{1}{2}}.
\end{equation*}%
So our calculus is not sharp: we get $\frac{1}{4}$ instead of $\frac{1}{2}.$

\bigskip

In the regularization Lemma \ref{RL} we have not assumed that $\sigma _{F}$
is invertible but we preferred to keep ${\mathbb{P}}(\det \sigma _{F}\leq
\eta ).$ We give now a variant under a strong non degeneracy assumption for $%
F:$ for every $p\geq 1$%
\begin{equation}
{\mathbb{E}}\big((\det \sigma _{F})^{-p}\big)\leq C_{p}<\infty .
\label{e7a}
\end{equation}%
denote%
\begin{equation}
\mathcal{Q}_{l}(F)=\mathcal{C}_{l,2}(F)({\mathbb{E}}(\det \sigma
_{F})^{-2l})^{1/2}  \label{e7d}
\end{equation}%
and $\mathcal{C}_{l}(F)$ given in (\ref{norm2'}).\ 

\begin{lemma}
\label{C2}Let $f\in C_{b}^{\infty}({\mathbb{R}}^{d})$ and $F\in\mathcal{E}%
^{d}$ such that (\ref{e7a}) holds. For every $q,m\in {\mathbb{N}}$\ there
exists a universal constant $C$ (depending on $q$ and $m$ only) such that
for every multi-index $\gamma$ with $\left\vert \gamma\right\vert \leq m,$
every $\delta>0$ and every $\eta>0$%
\begin{equation}
\left\vert {\mathbb{E}}(\partial_{\gamma}f(F))-{\mathbb{E}}(\partial_{\gamma
}f_{\delta}(F))\right\vert \leq C\delta^{q}\left\Vert f\right\Vert _{\infty }%
\mathcal{Q}_{q+m}(F).  \label{e7c}
\end{equation}
\end{lemma}

\textbf{Proof }We follow the same reasoning as in the previous proof and we
come back to (\ref{e7b}), but we do no more multiply with $\Psi _{\eta
}(\det \sigma _{F}):$ 
\begin{align*}
& {\mathbb{E}}(\partial _{\gamma }f(F))-{\mathbb{E}}(\partial _{\gamma
}f_{\delta }(F)) \\
& =\frac{1}{q!}\sum_{\left\vert \alpha \right\vert =q}\int_{0}^{1}\int
dz\phi _{\delta }(z){\mathbb{E}}\big(\partial _{\alpha }\partial _{\gamma
}f(F+\lambda z)\big)z^{\alpha }(1-\lambda )^{q}d\lambda .
\end{align*}%
Using the standard integration by parts formula (\ref{Iteration}) we obtain,
for some $p\geq 1$ 
\begin{equation*}
\left\vert {\mathbb{E}}(\partial _{\alpha }\partial _{\gamma }f(F+\lambda
z))\right\vert \leq ({\mathbb{E}}(\det \sigma _{F})^{-2(q+m)})^{1/2}\mathcal{%
C}_{q+m,2}(F)\left\Vert f\right\Vert _{\infty }.
\end{equation*}%
And by (\ref{e7e}) we conclude that%
\begin{equation*}
\left\vert {\mathbb{E}}(\partial _{\gamma }f(F))-{\mathbb{E}}(\partial
_{\gamma }f_{\delta }(F))\right\vert \leq \mathcal{Q}_{q+m}(F)\left\Vert
f\right\Vert _{\infty }\delta ^{q}.
\end{equation*}%
$\square $\medskip

In \cite{[CLT]} one gives the following more sophisticated version of the
the regularization lemma for smooth functions with polynomial growth. More
precisely we denote by $C_{p}^{\infty }({\mathbb{R}}^{d})$ the space of
smooth functions such that for every $q\in {\mathbb{N}}$ there exists $%
L_{q}(f)$ and $l_{q}(f)$ such that, for every multi index with $\left\vert
\alpha \right\vert \leq q$ and every $x\in {\mathbb{R}}^{d}$%
\begin{equation*}
\left\vert \partial _{\alpha }f(x)\right\vert \leq L_{q}(f)(1+\left\vert
x\right\vert )^{l_{q}(f)}.
\end{equation*}%
Then we have the following result (see \cite{[CLT]} Lemma 5.3)

\begin{lemma}
\label{RLbis} Let $F\in \mathcal{E}^{d}$ and $q,m\in {\mathbb{N}}.$ There
exists some constant $C\geq 1,$ depending on $d,m$ and $q$ only, such that
for every $f\in C_{\mbox{\rm{\scriptsize {pol}}}}^{q+m}({\mathbb{R}}^{d}),$
every multi-index $\gamma $ with $\left\vert \gamma \right\vert =m$ and
every $\eta ,\delta >0$ and $p>1$%
\begin{equation}
\begin{array}{l}
\displaystyle\left\vert {\mathbb{E}}(\partial _{\gamma }f(x+F))-{\mathbb{E}}%
(\partial _{\gamma }f_{\delta }(x+F))\right\vert \leq C(1+\left\Vert
F\right\Vert _{pl_{0}(f)})^{l_{0}(f)}(1+|x|)^{l_{m}(f)}\smallskip \\ 
\displaystyle\quad \times \Big(L_{m}(f)c_{l_{m}(f)}2^{l_{m}(f)}{\mathbb{P}}%
^{(p-1)/p}(\det \sigma _{F}\leq \eta )+2^{l_{0}(f)}c_{l_{0}(f)+q}\,L_{0}(f)%
\frac{\delta ^{q}}{\eta ^{2(q+m)}}\mathcal{C}_{q+m,1}(F)\Big).%
\end{array}
\label{e8}
\end{equation}
\end{lemma}

\subsection{Convergence in total variation distance}

Let us introduce the following distances:%
\begin{equation}
d_{k}(F,G)=\sup \{\left\vert {\mathbb{E}}(f(F))-{\mathbb{E}}%
(f(G))\right\vert :\sum_{\left\vert \alpha \right\vert =k}\left\Vert
\partial ^{\alpha }f\right\Vert _{\infty }\leq 1\}  \label{dk}
\end{equation}%
In the case $k=0$ this means that the test functions $f$ are just measurable
and bounded and in this case $d_{0}$ is the so called \textquotedblleft
total variation distance" that we will denote by $d_{TV}.$ Another
interesting distance is the \textquotedblleft Wasserstein distance" 
\begin{equation*}
d_{W}(F,G)=d_{1}(F,G)=\sup \{\left\vert {\mathbb{E}}(f(F))-{\mathbb{E}}%
(f(G))\right\vert :\left\Vert \nabla f\right\Vert _{\infty }\leq 1\}.
\end{equation*}%
In many problems the estimate of the error involves some Taylor type
extensions and then the test functions have to be differentiable and the
norms of the derivatives come on. So we are able to estimate $d_{k}$ for
some $k.$ And then one asks about the possibility to obtain estimates for
measurable test functions, as in total variation distance. And one may use
the regularization lemma presented before in order to do it. We give several
forms of such a result. But as we will show, the regularization lemma can be
applied to other distances. As an example, in the sequel we consider also
the following distance between random vectors in ${\mathbb{R}}^{d}$: we set 
\begin{equation}
d_{CF}(F,G)=\sup \{\big|{\mathbb{E}}(e^{\mathbf{i}\<\vartheta ,F\>})-{%
\mathbb{E}}(e^{\mathbf{i}\<\vartheta ,G\>})\big|\,:\,\vartheta \in {\mathbb{R%
}}^{d}\}.  \label{dCF}
\end{equation}%
So, $d_{CF}$ is the maximum distance between the characteristic functions of 
$F$ and $G$. There are many situations where it is easier to obtain bounds
on the difference of characteristic functions, especially when the targets
in consideration are infinitely divisible. One may for instance consult \cite%
{[AH]} for an introduction to Stein's method theory for this kind of
distribution. Again, the regularization lemma allows one to pass from such
distance to the distance in total variation.

The key remark which allows to use the regularization lemma is the following.

\begin{lemma}
\label{cf} Let $\phi _{\delta }$ be the super-kernel introduced in the
previous section and let $F,G$ denote random vectors in ${\mathbb{R}}^{d}$.
We also fix a multi-index $\beta $\ with $\left\vert \beta \right\vert =r$
(including the void multi-index, in which case $r=0).$ Then for every $k\in {%
\mathbb{N}}$ there exists a constant $C$ depending on $k$ only such that for
every $f\in C_{b}({\mathbb{R}}^{d})$ 
\begin{equation}
\left\vert {\mathbb{E}}(\partial ^{\beta }(f\ast \phi _{\delta })(F))-{%
\mathbb{E}}(\partial ^{\beta }(f\ast \phi _{\delta })(G))\right\vert \leq
C\delta ^{-(k+r)}d_{k}(F,G)\times \left\Vert f\right\Vert _{\infty }.
\label{R1}
\end{equation}%
We also have, for a universal constant $C>0$, 
\begin{equation}
\left\vert {\mathbb{E}}(\partial ^{\beta }(f\ast \phi _{\delta })(F))-{%
\mathbb{E}}(\partial ^{\beta }(f\ast \phi _{\delta })(G))\right\vert \leq
C\delta ^{-(2d+r)}d_{CF}(F,G)\times \left\Vert f\right\Vert _{1}.  \label{R2}
\end{equation}%
And moreover, for every $\varepsilon >0$ one may find $C$ (depending on $%
\varepsilon $) such that%
\begin{equation}
\left\vert {\mathbb{E}}(\partial ^{\beta }(f\ast \phi _{\delta })(F))-{%
\mathbb{E}}(\partial ^{\beta }(f\ast \phi _{\delta })(G))\right\vert \leq
C\delta ^{-(2d+r)}d_{CF}(F,G)^{1-\varepsilon }\times \left\Vert f\right\Vert
_{\infty }.  \label{R3}
\end{equation}
\end{lemma}

\textbf{Proof} (\ref{R1}) is an immediate consequence of the
definition of $d_{k}$ and of $\left\Vert \partial ^{\alpha }(f\ast \phi
_{\delta })\right\Vert _{\infty }\leq C\delta ^{-\left\vert \alpha
\right\vert }.$

Let us prove (\ref{R2}). Take first $f\in \mathcal{S}$ (the Schwartz space).
Let $\mathcal{F}f(\xi )=\int f(x)e^{-2\pi \mathbf{i}\<x,\xi \>}dx$ denote
the Fourier transform of $f$ and $\mathcal{F}_{F}(\xi )={\mathbb{E}}%
(e^{i\left\langle F,\xi \right\rangle })$ be the characteristic function of $%
F$. Then, using the inverse $\mathcal{F}^{-1}$ of $\mathcal{F}$ 
\begin{equation*}
{\mathbb{E}}(f(F))={\mathbb{E}}(\mathcal{FF}^{-1}f(F))={\mathbb{E}}\int 
\mathcal{F}^{-1}f(\xi )e^{-2\pi \mathbf{i}\left\langle F,\xi \right\rangle
}d\xi =\int \mathcal{F}^{-1}f(\xi )\mathcal{F}_{F}(-2\pi \xi )d\xi .
\end{equation*}%
Take $\alpha $ such that $\partial ^{\alpha }=\partial _{1}^{2}\cdots
\partial _{d}^{2}$ and notice that, by using integration by parts, one
obtains $\mathcal{F}^{-1}\partial ^{\alpha }f(\xi )=(2\pi \mathbf{i}%
)^{2d}\prod_{i=1}^{d}\xi _{i}^{2}\mathcal{F}^{-1}f(\xi )$ so that%
\begin{equation*}
{\mathbb{E}}(f(F))=\frac{1}{(2\pi \mathbf{i})^{2d}}\int \prod_{i=1}^{d}\xi
_{i}^{-2}\mathcal{F}^{-1}\partial ^{\alpha }f(\xi )\mathcal{F}_{F}(-2\pi 
\mathbf{i}\xi )d\xi .
\end{equation*}%
It follows that 
\begin{equation*}
\left\vert {\mathbb{E}}(f(F))-{\mathbb{E}}(f(G))\right\vert \leq \left\Vert 
\mathcal{F}_{F}-\mathcal{F}_{G}\right\Vert _{\infty }\left\Vert \mathcal{F}%
^{-1}\partial ^{\alpha }f\right\Vert _{\infty }=d_{CF}(F,G)\left\Vert 
\mathcal{F}^{-1}\partial ^{\alpha }f\right\Vert _{\infty }.
\end{equation*}%
We will use this formula with $f\ast \phi _{\delta }$ instead of $f.$ One has%
\begin{eqnarray*}
\left\Vert \mathcal{F}^{-1}\partial ^{\alpha }\partial ^{\beta }(f\ast \phi
_{\delta })\right\Vert _{\infty } &=&\delta ^{-(2d+r)}\left\Vert \mathcal{F}%
^{-1}(f\ast \partial ^{\alpha }\partial ^{\beta }\phi )\right\Vert _{\infty
}\leq \delta ^{-(2d+r)}\left\Vert f\ast \partial ^{\alpha }\partial ^{\beta
}\phi \right\Vert _{1} \\
&\leq &\delta ^{-(2d+r)}\left\Vert f\right\Vert _{1}\left\Vert \partial
^{\alpha }\partial ^{\beta }\phi \right\Vert _{1}\leq C\delta
^{-(2d+r)}\left\Vert f\right\Vert _{1}
\end{eqnarray*}%
so that, by inserting above, (\ref{R2}) is proved. In order to prove (\ref%
{R3}) we take a truncation function $\Psi _{M}\in C^{\infty }({\mathbb{R}}%
^{d})$ such that $1_{B_{M}(0)}\leq \Psi _{M}\leq 1_{B_{M+1}(0)}$ and we use (%
\ref{R2}) for $f\Psi _{M}.$ Since $\left\Vert f\Psi _{M}\right\Vert _{1}\leq
CM^{d}\left\Vert f\right\Vert _{\infty }$ we obtain 
\begin{equation*}
\left\vert {\mathbb{E}}(((f\Psi _{M})\ast \partial ^{\beta }\phi _{\delta
}(F))-{\mathbb{E}}((f\Psi _{M})\ast \partial ^{\beta }\phi _{\delta
}(G))\right\vert \leq C\delta ^{-(2d+r)}M^{d}\left\Vert f\right\Vert
_{\infty }d_{CF}(F,G).
\end{equation*}%
Notice that for $q$ one has%
\begin{equation*}
\int \left\vert y\right\vert ^{q}\left\vert \partial ^{\beta }\phi _{\delta
}(y)\right\vert dy=\delta ^{q-r}\int \left\vert y\right\vert ^{q}\left\vert
\partial ^{\beta }\phi (y)\right\vert dy.
\end{equation*}%
Moreover, since $F$ has finite moments of any order, one has, 
\begin{equation*}
{\mathbb{P}}(\left\vert F-y\right\vert \geq M+1)\leq M^{-q}{\mathbb{E}}%
(\left\vert F-y\right\vert ^{q})\leq CM^{-q}(1+\left\vert y\right\vert ^{q}).
\end{equation*}
So, for every $q\geq r$%
\begin{eqnarray*}
\left\vert {\mathbb{E}}(f(1-\Psi _{M})\ast \partial ^{\beta }\phi _{\delta
}(F))\right\vert &\leq &\left\Vert f\right\Vert _{\infty }\int {\mathbb{P}}%
(\left\vert F-y\right\vert \geq M+1)\left\vert \partial ^{\beta }\phi
_{\delta }(y)\right\vert dy \\
&\leq &C\delta ^{-r}\left\Vert f\right\Vert _{\infty }M^{-q}.
\end{eqnarray*}%
The same is true for $G$ and then, combining the two previous estimates we
obtain, for every $q\in {\mathbb{N}}$ 
\begin{equation*}
\left\vert {\mathbb{E}}(f\ast \partial ^{\beta }\phi _{\delta }(F))-{\mathbb{%
E}}(f\ast \partial ^{\beta }\phi _{\delta }(G))\right\vert \leq C\delta
^{-(2d+r)}\left\Vert f\right\Vert _{\infty }(d_{CF}(F,G)M^{d}+M^{-q}).
\end{equation*}%
We optimize over $M$ and we obtain (\ref{R3}). 
$\square $

\medskip

We first consider the case in which both $F$ and $G$ satisfy the weak non
degeneracy condition in (\ref{e8a}).

\begin{lemma}
\label{l1}Let $F,G\in \mathcal{E}^{d}$ be such that (\ref{e8a}) holds true
for some fixed $\kappa >0.$ Then for every $k\in {\mathbb{N}}$%
\begin{equation}
d_{TV}(F,G)\leq C\times (C_{\kappa ,q}(F)+C_{\kappa ,q}(G))^{\frac{k}{k+a}%
}d_{k}(F,G)^{\frac{a}{k+a}}  \label{e8c}
\end{equation}%
with%
\begin{equation}
C_{\kappa ,q}(F)=\theta _{\kappa }(F)^{\frac{2q}{\kappa +2q}}\mathcal{C}%
_{q,2}(F)^{\frac{\kappa }{\kappa +2q}},\qquad a=\frac{\kappa q}{\kappa +2q}.
\label{e8d}
\end{equation}%
Moreover, for every $\varepsilon >0$%
\begin{equation}
d_{TV}(F,G)\leq C\times (C_{\kappa ,q}(F)+C_{\kappa ,q}(G))^{\frac{2d}{2d+a}%
}d_{CF}(F,G)^{\frac{a(1-\varepsilon )}{2d+a}}  \label{e8c'}
\end{equation}
\end{lemma}

\textbf{Proof. }By (\ref{R1}) (with $r=0)$%
\begin{equation*}
\left\vert {\mathbb{E}}(f\ast \phi _{\delta }(F))-{\mathbb{E}}(f\ast \phi
_{\delta }(G))\right\vert \leq \left\Vert f\right\Vert _{\infty
}d_{k}(F,G)\delta ^{-k},
\end{equation*}%
so, using (\ref{e8b}) we get 
\begin{equation*}
\left\vert {\mathbb{E}}(f(F))-{\mathbb{E}}(f(G))\right\vert \leq C\left\Vert
f\right\Vert _{\infty }((C_{\kappa ,q}(F)+C_{\kappa ,q}(G))\delta ^{\frac{%
\kappa q}{\kappa +2q}}+d_{k}(F,G)\delta ^{-k}).
\end{equation*}%
We optimize over $\delta $ and we get (\ref{e8c}). The proof of (\ref{e8c'})
is the same, but one employes (\ref{R3}). $\square $\medskip

We give now a result under the strong non degeneracy condition both for $F$
and $G$: $(\det \sigma _{F})^{-1}, (\det \sigma _{G})^{-1}\in \cap _{p\geq
1}L^{p}.$

\begin{lemma}
\label{l2}Let $F,G\in \mathcal{E}$ be such that $\mathcal{Q}_{q}(F)+\mathcal{%
Q}_{q}(G)<\infty $ for every $q\in {\mathbb{N}}$ (see (\ref{e7d})). Then,
for every $q,k\in {\mathbb{N}}$ there exists a constant $C$ (depending on $q$
and $k$ only) such that 
\begin{equation}
d_{TV}(F,G)\leq C(\mathcal{Q}_{q}(F)+\mathcal{Q}_{q}(G))^{\frac{k}{q+k}%
}\times d_{k}^{\frac{q}{q+k}}(F,G).  \label{e12a}
\end{equation}%
In particular, for every $\varepsilon >0$ 
\begin{equation}
d_{TV}(F,G)\leq C(\mathcal{Q}_{q(\varepsilon )}(F)+\mathcal{Q}%
_{q(\varepsilon )}(G))^{\varepsilon }\times d_{k}^{1-\varepsilon }(F,G)
\label{e12b}
\end{equation}%
with $q(\varepsilon )=k(\frac{1}{\varepsilon }-1).$ Moreover, with $%
q(\varepsilon )=2d(\frac{1}{\varepsilon }-1),$ 
\begin{equation}
d_{TV}(F,G)\leq C(\mathcal{Q}_{q(\varepsilon )}(F)+\mathcal{Q}%
_{q(\varepsilon )}(G))^{\varepsilon }\times d_{CF}^{(1-\varepsilon
)^{2}}(F,G)  \label{e12b'}
\end{equation}
\end{lemma}

\textbf{Proof} Let $f\in C_{b}^{\infty }({\mathbb{R}})$ and $\delta >0.$
Using (\ref{e7c})%
\begin{equation*}
\left\vert {\mathbb{E}}(f(F))-{\mathbb{E}}(f\ast \phi _{\delta
}(F))\right\vert \leq C\delta ^{q}\left\Vert f\right\Vert _{\infty }\mathcal{%
Q}_{q}(F)
\end{equation*}%
and a similar estimate holds for $G.$ Moreover by (\ref{R1}) 
\begin{equation*}
\left\vert {\mathbb{E}}(f_{\delta }(F))-{\mathbb{E}}(f\ast \phi _{\delta
}(G))\right\vert \leq \frac{C}{\delta ^{k}}\left\Vert f\right\Vert _{\infty
}d_{k}(F,G)
\end{equation*}%
so that 
\begin{equation*}
\left\vert {\mathbb{E}}(f(F))-{\mathbb{E}}(f(G))\right\vert \leq C\left\Vert
f\right\Vert _{\infty }(\delta ^{q}(\mathcal{Q}_{q}(F)+\mathcal{Q}_{q}(F))+%
\frac{1}{\delta ^{k}}d_{k}(F,G)).
\end{equation*}%
We optimize on\ $\delta $ and we get (\ref{e12a}). Using (\ref{R3})$\ $we
obtain, in the same way, (\ref{e12b'}).\ $\square $

\begin{remark}
Compare with Corollary 2.8 pg 11/33 in \cite{[BCDi]} - there we have $%
d_{k}^{1/(k+1)}(F,F_{n}).$ So it is much less good there. The reason is that
we have now a much stronger regularization lemma. Compare the estimate
(2.29) pg 8/33 in \cite{[BCDi]} with (\ref{e7}) here: here we have
\textquotedblleft for every q" and this is what gives the much better result.
\end{remark}

We finish this section with a variant of the previous Lemma: now we assume
the strong non degeneracy condition $(\det\sigma_{F})^{-1}\in\cap_{p%
\geq1}L^{p}$ for $F$ but we assume no non degeneracy condition on $G.$ Then
we get the following:

\begin{proposition}
\label{l3}Let $F,G\in \mathcal{E}$ be such that $\mathcal{Q}_{q}(F)+\mathcal{%
C}_{q,1}(G)<\infty $ for every $q\in {\mathbb{N}}$ (see (\ref{e7d})). Then,
for every $p,p^{\prime }\in {\mathbb{N}}$ and every $\varepsilon >0$ there
exists a constant $C$ (depending on $\varepsilon ,p,p^{\prime })$ such that 
\begin{equation}
d_{TV}(F,G)\leq C\times C_{\varepsilon }(F,G)\times (d_{p}(F,G)+d_{p^{\prime
}}(\det \sigma _{F},\det \sigma _{G}))^{1-\varepsilon }  \label{e12c}
\end{equation}%
with 
\begin{equation}
C_{\varepsilon }(F,G)=1+\mathcal{Q}_{q(\varepsilon )}(F)+\mathcal{C}%
_{q(\varepsilon ),1}(G),\qquad q(\varepsilon )=\max \{[4p/\varepsilon
]+1,[p^{\prime }/2\varepsilon ]+1\}.  \label{e12d}
\end{equation}%
Moreover, with $q(\varepsilon )=\max \{[8d/\varepsilon ]+1,[p^{\prime
}/2\varepsilon ]+1\}$ one has%
\begin{equation}
d_{TV}(F,G)\leq C\times C_{\varepsilon }(F,G)\times
(d_{CF}(F,G)+d_{p^{\prime }}(\det \sigma _{F},\det \sigma
_{G}))^{1-\varepsilon }.  \label{e12c'}
\end{equation}
\end{proposition}

\textbf{Proof.} We denote 
\begin{equation*}
d_{p,p^{\prime}}=d_{p}(F,G)+d_{p^{\prime}}(\det\sigma_{F},\det\sigma_{G}).
\end{equation*}

Take $\eta >0$ and take $\Phi _{\eta }\in C_{b}^{\infty }({\mathbb{R}}_{+})$
such that $1_{(0,\eta )}\leq \Phi _{\eta }\leq 1_{(0,2\eta )}$ and $%
\Vert \Phi _{\eta }^{(k)}\Vert _{\infty }\leq C\eta ^{-k}.$ We
recall that $\sigma _{G}$ is the Malliavin covariance matrix of $G$ and we
write%
\begin{align}
{\mathbb{P}}(\det \sigma _{G}& \leq \eta )\leq {\mathbb{E}}(\Phi _{\eta
}(\det \sigma _{G}))\leq {\mathbb{E}}(\Phi _{\eta }(\det \sigma
_{F}))+\left\vert {\mathbb{E}}(\Phi _{\eta }(\det \sigma _{G}))-{\mathbb{E}}%
(\Phi _{\eta }(\det \sigma _{F}))\right\vert   \label{qui} \\
& \leq {\mathbb{P}}(\sigma _{F}\leq 2\eta )+C\eta ^{-p^{\prime
}}d_{p^{\prime }}(\det \sigma _{F},\det \sigma _{G})  \notag \\
& \leq C(\mathcal{Q}_{\rho }(F)\eta ^{\rho }+\eta ^{-p^{\prime
}}d_{p,p^{\prime }})  \notag
\end{align}%
the last inequality being true for every $\rho .$ Then, using the
regularization lemma, we get for every $\delta >0$ and every $q\in {\mathbb{N%
}}$%
\begin{align*}
\left\vert {\mathbb{E}}(f(G))-{\mathbb{E}}(f_{\delta }(G))\right\vert & \leq
C\left\Vert f\right\Vert _{\infty }\Big({\mathbb{P}}(\det \sigma _{G}\leq \eta )+%
\frac{\delta ^{q}}{\eta ^{2q}}\mathcal{C}_{q,1}(G)\Big) \\
& \leq C\left\Vert f\right\Vert _{\infty }(\mathcal{Q}_{\rho }(F)+\mathcal{C}%
_{q}(G))\Big(\eta ^{\rho }+\eta ^{-p^{\prime }}d_{p,p^{\prime }}+\frac{\delta
^{q}}{\eta ^{2q}}\Big).
\end{align*}%
We also have ${\mathbb{P}}(\sigma _{F}\leq \eta )\leq \mathcal{Q}_{\rho
}(F)\eta ^{\rho }$ for every $\rho $ so that the regularization lemma for $F$
gives%
\begin{equation*}
\left\vert {\mathbb{E}}(f(F))-{\mathbb{E}}(f_{\delta }(F))\right\vert \leq
C\left\Vert f\right\Vert _{\infty }\mathcal{Q}_{\rho }(F)(\eta ^{\rho }+%
\frac{\delta ^{q}}{\eta ^{2q}}).
\end{equation*}%
On the other hand%
\begin{equation*}
\left\vert {\mathbb{E}}(f_{\delta }(F))-{\mathbb{E}}(f_{\delta
}(G))\right\vert \leq C\left\Vert f\right\Vert _{\infty }\mathcal{\delta }%
^{-p}d_{p}(F,G).
\end{equation*}%
Putting these together%
\begin{equation*}
\left\vert {\mathbb{E}}(f(F))-{\mathbb{E}}(f(G))\right\vert \leq C\left\Vert
f\right\Vert _{\infty }(1+\mathcal{Q}_{\rho }(F)+\mathcal{C}_{q}(G))(%
\mathcal{\delta }^{-p}d_{p,p^{\prime }}+\eta ^{\rho }+\eta ^{-p^{\prime
}}d_{p,p^{\prime }}+\frac{\delta ^{q}}{\eta ^{2q}}).
\end{equation*}%
We fix now $\varepsilon >0$ and we choose $\delta $ and $\eta .$ First take $%
\delta =d_{p,p^{\prime }}^{2\varepsilon }$ so that 
\begin{equation*}
\mathcal{\delta }^{-p}d_{p,p^{\prime }}=d_{p,p^{\prime }}^{1-2p\varepsilon }.
\end{equation*}%
Take $\eta =d_{p,p^{\prime }}^{\varepsilon /2}$ so that $d_{p,p^{\prime
}}^{\varepsilon }\times \eta ^{-2}=1$ and consequently%
\begin{equation*}
\frac{\delta ^{q}}{\eta ^{2q}}=\frac{d_{p,p^{\prime }}^{2q\varepsilon }}{%
\eta ^{2q}}=d_{p,p^{\prime }}^{q\varepsilon }\times \frac{d_{p,p^{\prime
}}^{q\varepsilon }}{\eta ^{2q}}=d_{p,p^{\prime }}^{q\varepsilon }.
\end{equation*}%
We also have%
\begin{equation*}
\eta ^{-p^{\prime }}d_{p,p^{\prime }}=d_{p,p^{\prime }}^{1-p^{\prime
}\varepsilon /2}\qquad and\qquad \eta ^{\rho }=d_{p;p^{\prime }}^{\rho
\varepsilon /2}.
\end{equation*}%
Choose $q=1/\varepsilon $ and $\rho =2/\varepsilon $ in order to obtain 
\begin{equation*}
\eta ^{\rho }+\frac{\delta ^{q}}{\eta ^{2q}}\leq 2d_{p,p^{\prime }}.
\end{equation*}%
Then%
\begin{equation*}
\left\vert {\mathbb{E}}(f(F))-{\mathbb{E}}(f(G))\right\vert \leq C\left\Vert
f\right\Vert _{\infty }(1+\mathcal{Q}_{[2/\varepsilon ]+1}(F)+\mathcal{C}%
_{[\varepsilon ^{-1}]+1,1}(G))(d_{p,p^{\prime }}^{1-2p\varepsilon
}+d_{p,p^{\prime }}^{1-p^{\prime }\varepsilon /2}+d_{p,p^{\prime }}).
\end{equation*}%
Take now $\overline{\varepsilon }=\max \{2p\varepsilon ,p^{\prime
}\varepsilon /2\}$ and $q(\overline{\varepsilon })=\max \{[4p/\overline{%
\varepsilon }]+1,[p^{\prime }/2\overline{\varepsilon }]+1\}.$ The above
estimate reads%
\begin{equation*}
\left\vert {\mathbb{E}}(f(F))-{\mathbb{E}}(f(G))\right\vert \leq C\left\Vert
f\right\Vert _{\infty }(1+\mathcal{Q}_{q(\overline{\varepsilon })}(F)+%
\mathcal{C}_{q(\overline{\varepsilon }),1}(G))\times d_{p,p^{\prime }}^{1-%
\overline{\varepsilon }}
\end{equation*}%
so (\ref{e12c}) is proved. In order to prove (\ref{e12c'}) proceed as before
but we use (\ref{R3}) instead of (\ref{R1}). $\square $

\medskip In inequality (\ref{qui}) one can replace $\det \sigma_F$ with any
other random variable $H\geq 0$ such that ${\mathbb{P}}(H\leq \eta)\leq
c_\kappa \eta^\kappa$ for every $\kappa$ and $\eta>0$, that is, $H^{-1}\in
\cap_{p} L^p$. So, Proposition \ref{l3} can be reformulated as follows:

\begin{proposition}
\label{l3bis} Let $F,G\in \mathcal{E}$ be such that $\mathcal{Q}_{q}(F)+%
\mathcal{C}_{q,1}(G)<\infty $ for every $q\in {\mathbb{N}}$ (see (\ref{e7d}%
)). Let $H\geq 0$ be a r.v. such that $H^{-1}\in \cap _{p}L^{p}$. Then, for
every $p,p^{\prime }\in {\mathbb{N}}$ and every $\varepsilon >0$ there
exists a constant $C$ (depending on $\varepsilon ,p,p^{\prime })$ such that 
\begin{equation}
d_{TV}(F,G)\leq C\times C_{\varepsilon }(F,G,H)\times
(d_{p}(F,G)+d_{p^{\prime }}(\det \sigma _{F},H))^{1-\varepsilon },
\label{e12c}
\end{equation}%
with $C_{\varepsilon }(F,G,H)=C_{\varepsilon }(F,G)+\Vert
H^{-1}\Vert _{2/\varepsilon }$\ with $C_{\varepsilon }(F,G)$ given in (%
\ref{e12d}). Moreover, 
\begin{equation}
d_{TV}(F,G)\leq C\times C_{\varepsilon }(F,G,H)\times
(d_{CF}(F,G)+d_{p^{\prime }}(\det \sigma _{F},H))^{1-\varepsilon }.
\label{e12c'bis}
\end{equation}
\end{proposition}

\begin{remark}
Proposition \ref{l3} essentially says that if $(F_{n},\det\sigma_{F_{n}})%
\rightarrow(F,\det\sigma_{F})$ in some ``smooth distance" $d_{p}$ (for
example in the Wasserstein distance $d_{W}$) then $F_{n}\rightarrow F$ in
total variation distance. And Proposition \ref{l3bis} says that it is not
necessary that $\sigma_{F_{n}}$ converges to $\sigma_F$: one can also have $%
(F_{n},\det\sigma_{F_{n}})\rightarrow(F,H)$. And one obtains the estimate of
the speed of convergence: one looses a little bit because we have the power $%
1-\varepsilon.$ The striking fact is the we do not need the non degeneracy
condition for $F_{n}$ but only for $F.$
\end{remark}

\begin{remark}
For the use of Proposition \ref{l3}, in concrete applications it may be
difficult to compute $d_{p^{\prime }}(\det \sigma _{F},\det \sigma _{G}).$
Then we are obliged to come back to \textquotedblleft strong distances": we
have $d_{1}(\det \sigma _{F},\det \sigma _{G})\leq (\left\Vert DF\right\Vert
^{d-1}+\left\Vert DG\right\Vert ^{d-1})\left\Vert DF-DG\right\Vert $ so we
take $p^{\prime }=1$ and we get 
\begin{equation}
d_{TV}(F,G)\leq C\times C_{\varepsilon }(F,G)\times (d_{p}(F,G)+\left\Vert
DF-DG\right\Vert )^{1-\varepsilon }  \label{e12e}
\end{equation}%
Notice however that here we loose the right order of convergence. For
example, in the case of the convergence of the Euler scheme $X_{t}^{n}$ to
the diffusion process $X_{t}$ we have $d_{1}(X_{t}^{n},X_{t})\leq \frac{C}{n}
$ but $\left\Vert DX_{t}^{n}-DX_{t}\right\Vert \sim \frac{C}{n^{1/2}}.$ This
is because we deal with the weak convergence in the first case and with the
strong convergence in the second one. So we pass from $\frac{1}{n}$ to $%
\frac{1}{n^{1/2}}.$ If we have an ellipticity property, then we do no need
to estimate $\left\Vert DX_{t}^{n}-DX_{t}\right\Vert $ so we are at level $%
\frac{1}{n}.$
\end{remark}

\subsection{The distance between density functions}

We give here an immediate application of the regularization lemma under the
strong non degeneracy condition (Lemma \ref{C2}): we estimate the distance
between density functions.

\begin{proposition}
\label{P}Let $F,G\in \mathcal{E}^{d}$ be such that $\mathcal{Q}_{q}(F)+ 
\mathcal{Q}_{q}(G)<\infty $ for every $q\in {\mathbb{N}}$ (see (\ref{e7d})).
Then for every $k\in N,$ every multi index $\alpha =(\alpha _{1},...,\alpha 
_{m})$ and every $\varepsilon >0$ there exists some constants $C$ and $q$ 
(depending on $\varepsilon ,k$ and on $m)$ such that  
\begin{equation}
\left\vert {\mathbb{E}}(\partial ^{\alpha }f(F))-{\mathbb{E}}(\partial
^{\alpha }f(G))\right\vert \leq Cd_{k}^{1-\varepsilon }(F,G)\left\Vert
f\right\Vert _{\infty }(\mathcal{Q}_{q}(F)+\mathcal{Q}_{q}(G))^{\varepsilon }
\label{e13}
\end{equation}
In particular if $p_{F}$ respectively $p_{G}$ are the density functions for $
F$ respectively $G$ then, for every $x\in {\mathbb{R}}^{d}$  
\begin{equation}
\left\vert \partial ^{\alpha }p_{F}(x)-\partial ^{\alpha
}p_{G}(x))\right\vert \leq Cd_{k}^{1-\varepsilon }(F,G)(\mathcal{Q}_{q}(F)+ 
\mathcal{Q}_{q}(G))^{\varepsilon }  \label{e14}
\end{equation}
The same estimates hold with $d_{k}^{1-\varepsilon }(F,G)$ replaced by $
d_{CF}^{1-\varepsilon }(F,G).$
\end{proposition}

\begin{remark}
Compare with the estimate (2.53) pg 14/33 in \cite{[BCDi]} : here the 
estimate is much better because, using $d_{1}$ for example, we have just $
d_{1}^{1-\varepsilon }(F,G)\leq {\mathbb{E}}(\left\vert F-G\right\vert 
))^{1-\varepsilon }$ and the Sobolev norms are not involved (as it is the 
case in \cite{[BCDi]}).
\end{remark}

\textbf{Proof }We will prove just (\ref{e13}) because (\ref{e14}) follows by
standard regularization methods. To begin we use (\ref{e7c}) and we get%
\begin{equation*}
\left\vert {\mathbb{E}}(\partial ^{\alpha }f(F))-{\mathbb{E}}(\partial
^{\alpha }(f\ast \phi _{\delta })(F))\right\vert \leq C\delta ^{q}\left\Vert
f\right\Vert _{\infty }\mathcal{Q}_{q+m}(F)
\end{equation*}%
and a similar estimate for $G.$ Moreover, using (\ref{R1})

\begin{equation*}
\left\vert {\mathbb{E}}(\partial ^{\alpha }(f\ast \phi _{\delta })(F))-{%
\mathbb{E}}(\partial ^{\alpha }(f\ast \phi _{\delta })(G))\right\vert \leq
C\delta ^{-(k+m)}d_{k}(F,G)\times \left\Vert f\right\Vert _{\infty }
\end{equation*}%
so that 
\begin{equation*}
\left\vert {\mathbb{E}}(\partial ^{\alpha }f(F))-{\mathbb{E}}(\partial
^{\alpha }f(G))\right\vert \leq C(\delta ^{-(k+m)}d_{k}(F,G)+\delta ^{q}(%
\mathcal{Q}_{q+m}(F)+\mathcal{Q}_{q+m}(G)))\left\Vert f\right\Vert _{\infty
}.
\end{equation*}%
Optimizing over $\delta $ we get (\ref{e13}). $\square $

\bigskip 

\section{Examples}

\subsection{Euler Scheme}

In this section we discuss the convergence in total variation distance of
the Euler scheme.

We consider the $d$ dimensional diffusion process 
\begin{equation*}
X_{t}=x+\sum_{j=1}^{m}\int_{0}^{t}\sigma
_{j}(X_{s})dB_{s}^{j}+\int_{0}^{t}b(X_{s})ds
\end{equation*}%
where $\sigma _{j}\in C_{b}^{\infty }({\mathbb{R}}^{d},{\mathbb{R}}%
^{d}),j=1,...,m$ and $b\in C_{b}^{\infty }({\mathbb{R}}^{d},{\mathbb{R}}%
^{d}).$ We are concerned with the Euler scheme defined in the following way.
We fix $n\in {\mathbb{N}},$ we define $\tau (t)=\frac{k}{n}$ for $\frac{k}{n}%
\leq t<\frac{k+1}{n}$ and then the Euler scheme is given by 
\begin{equation*}
X_{t}^{n}=x+\sum_{j=1}^{m}\int_{0}^{t}\sigma _{j}(X_{\tau
(s)}^{n})dB_{s}^{j}+\int_{0}^{t}b(X_{\tau (s)}^{n})ds.
\end{equation*}%
In this framework there are two types of errors which are of interest:
first, the \textquotedblleft strong error" is $\left\Vert
X_{t}-X_{t}^{n}\right\Vert _{p}\leq \frac{C}{\sqrt{n}}.$ And moreover, the
weak error: for $f\in C_{b}^{6}({\mathbb{R}}^{d}),$ Talay and Tubaro proved
in \cite{[TT]} that%
\begin{equation}
\left\vert {\mathbb{E}}(f(X_{t}))-{\mathbb{E}}(f(X_{t}^{n}))\right\vert \leq
C\left\Vert f\right\Vert _{6,\infty }\times \frac{1}{n}.  \label{es1}
\end{equation}%
Then one is interested to obtain the above estimate for measurable and
bounded functions, so to replace $\left\Vert f\right\Vert _{6,\infty }$ by $%
\left\Vert f\right\Vert _{\infty }$ in the above estimate (this means to
obtain the estimate of the error in total variation distance). This has been
done by Bally and Talay in \cite{[BT]} and by Guyon in \cite{[G]} by using
Malliavin calculus (another approach, based on the parametrix method has
been given by Konackov and Memen \cite{[KM]}, in the elliptic case). In
order to do this one has to assume a non degeneracy hypothesis. We construct
the Lie algebra associated to the coefficients of the above $SDE:$%
\begin{align*}
\mathcal{A}_{0}& =\{\sigma _{1},...,\sigma _{m}\}, \\
\mathcal{A}_{k}& =\{[\sigma _{1},\psi ],...,[\sigma _{m},\psi ],[b,\psi
]:\psi \in \mathcal{A}_{k-1}\}
\end{align*}%
where $[\phi ,\psi ]=\left\langle \phi ,\nabla \psi \right\rangle
-\left\langle \psi ,\nabla \phi \right\rangle $ is the Lie bracket We also
denote $\mathcal{A}_{k}(x)=\{\phi (x),\phi \in \mathcal{A}_{k}\}.$ Then we
have two types of non degeneracy conditions: the ellipticity condition in $x$
means that $\sigma _{1}(x),...,\sigma _{m}(x)$ span ${\mathbb{R}}^{d}.$ This
is also equivalent with the fact that $\sigma \sigma ^{\ast }(x)$ is
invertible The second condition, much less strong, is that $\cup _{k\in {%
\mathbb{N}}}\mathcal{A}_{k}(x)$ span ${\mathbb{R}}^{d}.$ This is the so
called H\"{o}rmander's condition. And H\"{o}rmander's theorem (proved by
Malliavin by a probabilistic approach) say that under this condition the law
of $X_{t}(x)$ is absolutely continuous and has a smooth density.

Let us come back to the estimate of the weak error in total variation
distance. J. Guyon proved in \cite{[G]} that if the uniform ellipticity
condition holds that is $\sigma \sigma ^{\ast }(x)\geq \lambda >0$ for every 
$x,$ then 
\begin{equation}
\left\vert {\mathbb{E}}(f(X_{t}(x)))-{\mathbb{E}}(f(X_{t}^{n}(x)))\right%
\vert \leq C\left\Vert f\right\Vert _{\infty }\times \frac{1}{n}.
\label{es3}
\end{equation}%
The estimate of the weak error in total variation distance, under the H\"{o}%
rmander condition, has been done in \cite{[BT]} under the \textquotedblleft
uniform H\"{o}rmander condition": there is a $k\in {\mathbb{N}}$ and some $%
\lambda >0$ such that 
\begin{equation*}
\Lambda (x):=\inf_{\left\vert \xi \right\vert =1}\sum_{\psi \in \mathcal{A}%
_{k}}\left\langle \psi (x),\xi \right\rangle ^{2}\geq \lambda .
\end{equation*}

But here a supplementary difficulty appears: the H\"{o}rmander assumption is
not sufficient in order to guarantee that Malliavin covariance matrix $%
\sigma_{X_{t}^{n}(x)}$ of the Euler scheme $X_{t}^{n}(x)$ is invertible (and
this was a crucial ingredient in the proof). In \cite{[BT]} this difficulty
has been bypassed by replacing $X_{t}^{n}(x)$ by the ``regularized version" $%
\overline{X}_{t}^{n}(x)=X_{t}^{n}(x)+\varepsilon_{n}\Delta$ where $\Delta$
is a standard normal random variable independent of the Brownian motion $B.$
From a simulation point of view this is really not a problem because this
just means to simulate one more random variable $\Delta.$ And one proves that%
\begin{equation}
\left\vert {\mathbb{E}}(f(X_{t}(x)))-{\mathbb{E}}(f(\overline{X}%
_{t}^{n}(x)))\right\vert \leq C\left\Vert f\right\Vert _{\infty}\times\frac{1%
}{n}.  \label{es2}
\end{equation}
Although from a practical point of view this has not big interest, the
following theoretical question remained open: is it possible to prove (\ref%
{es2}) for the real Euler scheme $X_{t}^{n}(x)$ (without the regularization
factor $\varepsilon_{n}\Delta)$ under the H\"{o}rmander condition?

Let us see what we may obtain using the results from the previous sections.
We use the standard Malliavin calculus so now $\mathcal{E}={\mathbb{D}}%
^{\infty}$ (see the notation in Nualart \cite{bib:[N]}). And under our
assumptions on the coefficients ($\sigma_{j},b\in C_{b}^{\infty}({\mathbb{R}}%
^{d},{\mathbb{R}}^{d})$) standard estimates yield, for every $q\in {\mathbb{N%
}}$%
\begin{equation*}
\mathcal{C}_{q}(X_{t}(x))+\sup_{n}\mathcal{C}_{q}(X_{t}^{n}(x))=C_{q}<\infty.
\end{equation*}
Moreover, if the ellipticity condition holds, one proves that $\det
\sigma_{X_{t}(x)}\geq\lambda(x)>0$ and $\sigma_{X_{t}^{n}(x)}\geq%
\lambda(x)>0 $ with $\lambda(x)$ independent of $n.$ It follows that 
\begin{equation*}
\mathcal{Q}_{q}(X_{t}(x))+\sup_{n}\mathcal{Q}_{q}(X_{t}^{n}(x))=Q_{q}<\infty.
\end{equation*}
So, using (\ref{e12b}) first and (\ref{es1}) then, we obtain for every $%
\varepsilon>0$ 
\begin{equation}
d_{TV}(X_{t},X_{t}^{n})\leq C\times
d_{6}^{1-\varepsilon}(X_{t},X_{t}^{n})\leq\frac{C}{n^{1-\varepsilon}}.
\label{es4}
\end{equation}
Comparing with the result of Guyon (\ref{es3}) we see that we have lost a
little bit because we have the power $1-\varepsilon$ instead of $1.$ This is
a structural drawback of our method which is based on optimization. However
there is a slight gain because we need just ellipticity in the starting
point $x$ and not uniform ellipticity.

Let us see now what we are able to say under H\"{o}rmander's condition. We
stress that we do not need the uniform H\"{o}rmander condition but only the
condition in the starting point $x:$ we just assume that $Span\{\cup _{k\in {%
\mathbb{N}}}\mathcal{A}_{k}(x)={\mathbb{R}}^{d}\}.$\ This is sufficient in
order to guarantee that $\det \sigma _{X_{t}(x)}>0$ and this is all we need.
As we mentioned above, we are no more able to prove that $\sigma
_{X_{t}^{n}(x)}\geq \lambda (x)>0$ so we have to use (\ref{e12c}) (together
with (\ref{es1})) : 
\begin{align*}
d_{TV}(X_{t},X_{t}^{n})& \leq C\times (d_{6}(X_{t},X_{t}^{n})+d_{p^{\prime
}}(\det \sigma _{X_{t}},\det \sigma _{X_{t}^{n}}))^{1-\varepsilon } \\
& \leq C\times \Big(\frac{1}{n}+d_{p}(\det \sigma _{X_{t}},\det \sigma
_{X_{t}^{n}})\Big)^{1-\varepsilon }.
\end{align*}%
Now we have to estimate $d_{p}(\det \sigma _{X_{t}},\det \sigma
_{X_{t}^{n}}).$ If we would be able to prove that for some $p\in {\mathbb{N}}
$ one has $d_{p}(\det \sigma _{X_{t}},\det \sigma _{X_{t}^{n}})\leq \frac{C}{%
n}$ then we come back to the same estimate as in the elliptic case. At a
first glance this seems reasonable, but taking things seriously this is not
so clear - we give up to answer this question here, and we just notice that
easy standard arguments give $\left\Vert \det \sigma _{X_{t}}-\det \sigma
_{X_{t}^{n}})\right\Vert _{1}\leq \frac{C}{\sqrt{n}}$ which yields $%
d_{1}(\det \sigma _{X_{t}},\det \sigma _{X_{t}^{n}})\leq \frac{C}{n^{1/2}}.$
Finally we obtain:%
\begin{equation}
d_{TV}(X_{t},X_{t}^{n})\leq \frac{C}{n^{\frac{1}{2}-\varepsilon }}.
\label{es5}
\end{equation}%
We conclude that we are still able to prove that $%
\lim_{n}d_{TV}(X_{t}(x),X_{t}^{n}(x))=0$ but we loose much on the speed of
convergence.

\subsection{Central Limit Theorem for Wiener chaoses}

Let us fix $1\leq q_1\leq q_2\cdots\le q_d$ a sequence of $d$ positive
integers. Let us consider here $F_n=(F_{1,n},\cdots,F_{d,n})$ a sequence of
random vectors such that for all $i\in\{1,\cdots,d\}$ and all $n\ge 1$, the
random variable $F_{i,n}$ belongs to the $q_i$--th Wiener chaos. We will
further assume that the covariance matrix of $F_n$ is the identity matrix
for every $n\ge 1$. A central result of Nourdin-Peccati theory established
in \cite{[NPS]} , provides an explicit bound in total variation between the
distribution $F_n$ and the distribution of a standard Gaussian vector, say $%
N=(N_1,\cdots,N_d)$:

\begin{equation}  \label{TV-multi-chaos}
d_{TV}\left(F,N\right)\le C \Phi\left({\mathbb{E}}\left(|F_n|^4\right)-%
\mathbb{E}\left(|N|^4\right)\right),\,\,\Phi(x)=|\log(x)|\sqrt{x}.
\end{equation}
Let us mention that an entropic result is actually proved in \cite{[NPS]}
and the previous bound is the corresponding total variation estimate which
is derived from Pinsker inequality.

\medskip

The proof of (\ref{TV-multi-chaos}) uses clever arguments from information
theory which are nevertheless rather specific to Gaussian targets. Our goal
here is to apply Proposition \ref{l3bis} to this situation and to compare
the bounds. First, from \cite{[NPR]} one has the following result regarding
the Wasserstein distance: 
\begin{equation}  \label{Wass-multi-chaos}
d_W\left(F,N\right)\le C \left(\sum_{i,j=1}^d {\mathbb{E}}%
\left(\left(\delta_{i,j}-\frac{1}{q_i}\Gamma[F_{i,n},F_{j,n}]%
\right)^2\right)\right)^{\frac 1 2},
\end{equation}
while from \cite{[NR]} one gets the bound 
\begin{equation}  \label{bound-fmt-multi}
\sum_{i,j=1}^d {\mathbb{E}}\left(\left(\delta_{i,j}-\frac{1}{q_i}\Gamma[%
F_{i,n},F_{j,n}]\right)^2\right) \le \mathbb{E}\left(|F_n|^4\right)-\mathbb{E%
}\left(|N|^4\right).
\end{equation}
Finally it is obvious that for some constant $C$ only depending on $d$ we
get 
\begin{equation}  \label{dist-Malliavin-matrix}
d_{W}\left(\det\sigma_{F_n},\det\mathbf{D}\right)\le C \left(\sum_{i,j=1}^d {%
\mathbb{E}}\left(\left(\delta_{i,j}-\frac{1}{q_i}\Gamma[F_{i,n},F_{j,n}]%
\right)^2\right)\right)^{\frac 1 2},
\end{equation}
with $\mathbf{D}=\text{Diag}(q_1,q_2,\cdots,q_d)$. Since $\mathbf{D}$ is a
deterministic invertible matrix and $F_n$ is uniformly bounded in $\mathbb{D}%
^\infty$, one can apply Proposition \ref{l3bis} with $H=\det \mathbf{D}$.
For every $\epsilon>0$, gathering the bounds (\ref{Wass-multi-chaos}), then (%
\ref{bound-fmt-multi}) and (\ref{dist-Malliavin-matrix}) lead to

\begin{equation}  \label{reg-lemma-chaosmulti}
d_{TV}\left(F,N\right)\le C_\epsilon \left(\mathbb{E}\left(|F_n|^4\right)-%
\mathbb{E}\left(|N|^4\right)\right)^{\frac{1}{2}-\epsilon}
\end{equation}
which, up to an arbitrarily small loss, retrieves the correct order of
magnitude in the total variation estimate of the so-called fourth moment
Theorem.

\appendix


\bigskip

\end{document}